\documentclass[11pt]{amsart}
\usepackage{amsfonts}
\usepackage{ifthen}
\usepackage{amsthm}
\usepackage{amsmath}
\usepackage{graphicx}
\usepackage{amscd,amssymb,amsthm}
\usepackage{graphicx}
\usepackage{epstopdf}
\usepackage{hyperref}

\newcounter{minutes}
\divide\time by 60
\newcounter{hours}
\multiply\time by 60 \addtocounter{minutes}{-\time}

\newtheorem{theorem}{Theorem}
\newtheorem{lemma}{Lemma}

\newcommand{\lp}{\mathcal{LP}}
\newcommand{\lpi}{\mathcal{LP}I}
\newcommand{\g}{\gamma}

\title{Zeros of some special entire functions}

\author[\'A. Baricz]{\'Arp\'ad Baricz}

\address{Department of Economics, Babe\c{s}-Bolyai University, Cluj-Napoca, Romania}

\address{Institute of Applied Mathematics, \'Obuda University, Budapest, Hungary}

\email{bariczocsi@yahoo.com}

\author[S. Singh]{Sanjeev Singh}
\address{Indian Statistical Institute, Chennai Centre, Chennai, India}
\email{sanjeevsinghiitm@gmail.com}



\begin{document}

\begin{abstract}
The real and complex zeros of some special entire functions such
as Wright, hyper-Bessel, and a special case of generalized hypergeometric functions are studied by using some
classical results of Laguerre, Obreschkhoff, P\'olya and Runckel. The obtained results extend the known theorem
 of Hurwitz on exact number of nonreal zeros of Bessel functions of the first kind.
Moreover, results on zeros of derivatives of Bessel functions and cross-product of Bessel functions
are also given, which are related to some recent open problems.
\end{abstract}

\def\thefootnote{}
\footnotetext{ \texttt{File:~\jobname .tex,
          printed: \number\year-0\number\month-\number\day,
          \thehours.\ifnum\theminutes<10{0}\fi\theminutes}
} \makeatletter\def\thefootnote{\@arabic\c@footnote}\makeatother

\maketitle

\section{Introduction and Main Results}
\setcounter{equation}{0}

Because of various applications in applied mathematics the zeros of Bessel functions of the first kind $J_{\nu}$ of order $\nu$ have been studied frequently. Lommel proved the reality of the zeros of Bessel functions for $\nu>-1$, while Hurwitz \cite{hur} completed the picture of the behavior of the zeros of Bessel functions of the first kind by determining the exact number of nonreal zeros for $\nu<-1.$ Many other interesting proof of Hurwitz's theorem were found; see the papers of Hilb \cite{hilb}, Obreschkoff \cite{obr}, P\'olya \cite{polya}, Hille and Szeg\H o \cite{hs}, Peyerimhoff \cite{peye}, Runckel \cite{runckel}, Ki and Kim \cite{kim}.  Hurwitz's original proof (based on Lommel polynomials) is quite long and difficult to read and Watson \cite{wa} even corrected some gaps in the proof. In this paper our aim is to point out that some results of Laguerre \cite{ti}, Obreschkoff \cite{obr}, P\'olya \cite{polya} and Runckel \cite{runckel} are useful to study the real and complex zeros of those special entire functions whose coefficients are expressed in terms of the reciprocal gamma function.

\subsection{Some classical results on zeros of entire functions}
A real entire function $q$ belongs to the Laguerre-P\'{o}lya class $\lp$ if it can be represented in the form
$$
q(z) = c z^{m} e^{-a z^{2} + \beta z} \prod_{k\geq1}
\left(1+\frac{z}{z_{k}}\right) e^{-\frac{z}{z_{k}}},
$$
with $c,$ $\beta,$ $z_{k} \in \mathbb{R},$ $a \geq 0,$ $m\in
\mathbb{N}_0,$ $\sum\limits_{k\geq1} z_{k}^{-2} < \infty.$ Here
$\mathbb{N}_0$ is the set of non-negative integers. Similarly, $w$  is  said to be of
type I in the Laguerre-P\'{o}lya class, written $w \in \lpi$,
if $w(z)$ or $w(-z)$ can be represented as
$$
w(z) = c z^{m} e^{\sigma z} \prod_{k\geq1}\left(1+\frac{z}{z_{k}}\right),
$$
with $c \in \mathbb{R},$ $\sigma \geq 0,$ $m \in
\mathbb{N}_0,$ $z_{k}>0,$ $\sum\limits_{k\geq1}z_{k}^{-1} < \infty.$
The class $\lp$ is the complement of the space of
polynomials whose all zeros are real in the
topology induced by the uniform convergence
on the compact sets of the complex plane, while $\lpi$ is the complement
of polynomials whose zeros are all real and posses a preassigned constant sign.
Given an entire function $\varphi$ in the form
$$\varphi(z) = \sum_{k\geq0}\gamma_{k} \frac{z^{k}}{k!},$$
its Jensen polynomial is defined by
$$
g_n(\varphi;z) = \sum_{j=0}^{n} {n\choose j} \g_{j}
z^j.
$$
Jensen proved the following result in \cite{Jen12}: The function $\varphi$ belongs to $\lp$ ($\lpi$) if and only if
the polynomials $g_n(\varphi;z)$ have only real zeros (real zeros of equal sign). Moreover, if for the function $\varphi\in\mathcal{LPI}$ we have $\gamma_k\geq0$ for all $k\in\mathbb{N}_0,$ we say that $\varphi\in\mathcal{LP}^{+}.$ Further information about the Laguerre-P\'olya class can be found in \cite{craven}, \cite{Obr} and \cite{DC}.

The following particular case of Laguerre' theorem contains sufficient condition on a special power series to have only real negative zeros. This result was used by P\'olya \cite{po} and it motivated P\'olya in writing \cite{polya}. Lemma \ref{lemlag} can be found in \cite[p. 270]{ti}, \cite[p. 186]{po} and also in \cite[p. 39]{dimi}.

\begin{lemma}[Laguerre, 1898]\label{lemlag}
If $\varphi$ is an entire function of order less than $2$ which takes real values
along the real axis and possesses only real negative zeros, then the
entire function $\sum\limits_{n\geq0}\frac{\varphi(n)}{n!}z^n$ has also real and negative zeros.
\end{lemma}

In order to shorten the proof of Hurwitz theorem on zeros of Bessel functions of the first kind, Obreschkoff \cite{obr} deduced the following result, which seems to be useful in proving the reality of zeros of special functions.

\begin{lemma}[Obreschkoff, 1929]\label{lemobr}
Let $q$ be an entire function of growth order $0$ or $1$, which has only real zeros, and has $s$ positive zeros. Then $\sum\limits_{n\geq0}\frac{(-1)^nq(2n)}{n!}z^{2n}$ has at most $2s$ complex zeros.
\end{lemma}

The following beautiful result of P\'olya \cite{polya} is useful in determining the exact number of nonreal zeros of some special entire functions. We note that the terminology of the original result of P\'olya \cite[p. 162]{polya} is a little bit different than our exposition: P\'olya uses the terminology of real oriented (reell gerichtet) and positively oriented (positiv gerichtet) functions, which are limits of polynomials (with real coefficients) having only real, and only positive real roots, respectively. In today's terminology these are members of $\lp$ and $\lpi.$

\begin{lemma}[P\'olya, 1929]\label{lempolya}
If the function $\sum\limits_{n\geq0}a_nz^n\in\mathcal{LPI}$ has nonzero roots and the function $G\in\mathcal{LP}$ has exactly $s$ simple roots in $[0,\infty)$ with the property that the distance between two arbitrary consecutive roots is not less than $1,$ then the function $\sum\limits_{n\geq0}a_nG(n)z^n$ has exactly $s$ non-positive roots.
\end{lemma}

The result of Runckel \cite[Theorem 4]{runckel}, stated in Lemma \ref{lemrun}, is also useful in proving the reality of zeros of some special functions.

\begin{lemma}[Runckel, 1969]\label{lemrun}
If $f(z)=\sum\limits_{n\geq0}a_nz^n$ can be represented as $f(z)=e^{az^2}h(z),$ where $a\leq0$ and $h$ is of type
$$h(z)=ce^{bz}\prod_{n\geq1}\left(1-\frac{z}{c_n}\right)e^{\frac{z}{c_n}},\ \ \ c,b\in\mathbb{R},\ \sum\limits_{n\geq1}|c_n|^{-2}<\infty;$$
$f$ has real zeros only (or no zeros at all), and $G$ is of type
\begin{equation}\label{typeA}G(z)=e^{\beta z}\prod_{n\geq1}\left(1+\frac{z}{\alpha_n}\right)e^{-\frac{z}{\alpha_n}},\ \ \ \alpha_n>0,\ \beta\in\mathbb{R},\ \sum\limits_{n\geq1}\alpha_n^{-2}<\infty,\end{equation} then the function $\sum\limits_{n\geq0}a_nG(n)z^n$ has real zeros only.
\end{lemma}

\subsection{Real and complex zeros of some special entire functions}

By using the above classical results of Laguerre, Obreschkhoff, P\'olya and Runckel our aim is to present some results related to real and complex zeros of some special entire functions such
as Wright, hyper-Bessel, a special case of generalized hypergeometric function, derivatives of Bessel function, product and cross-product of
Bessel and modified Bessel functions of the first kind. Moreover, we prove a result on a special function related to an open problem in \cite[Problem 6]{csord}. The results on Wright and hyper-Bessel functions extend
naturally the theorem of Hurwitz on zeros of Bessel functions of the first kind, which states that if $\nu\geq0,$ then
$$z\mapsto z^{\frac{\nu}{2}}J_{-\nu}(2\sqrt{z})=\sum_{n\geq0}\frac{(-1)^nz^n}{n!\Gamma(n-\nu+1)}$$
has exactly $[\nu]$ non-positive zeros.

Our first main result is a natural extension of Hurwitz's theorem on Bessel functions and it is about the zeros of the Wright function
$$\phi(\rho,\beta,z)=\sum_{n\geq0}\frac{z^n}{n!\Gamma(\rho n+\beta)},$$
where $\rho>-1$ and $\beta\in\mathbb{R}.$ We note that the special case of the first affirmation
in Theorem \ref{th1} for $\beta=1$ has been considered by Craven and Csordas \cite[Example 2.7]{craven}.

\begin{theorem}\label{th1}
If $\rho>0$ and $\beta>0,$ then all zeros of $\phi(\rho,\beta,-z)$ are real and positive. Moreover, if $0<\rho\leq 1$ and $\beta>0,$ then $\phi(\rho,-\beta,-z)$ has $[\beta]+1$ non-positive zeros.
\end{theorem}

The hyper-Bessel function (or a multi-index analogue of the Bessel function) is defined by
$$J_{\mathbf{\alpha_d}}(z)=\frac{\left(\frac{z}{d+1}\right)^{\alpha_1+\cdots+\alpha_d}}{\Gamma(\alpha_1+1)\cdots\Gamma(\alpha_d+1)}
{}_0F_d\left(-,{\bf \alpha_d+1};-\left(\frac{z}{d+1}\right)^{d+1}\right),$$
where $\mathbf{\alpha_d}=(\alpha_1,\dots,\alpha_d),$ $d\in\mathbb{N},$ and
$${}_pF_q({\bf a}; {\bf b}; z)=\sum_{n\geq0}\frac{(a_1)_n\cdots(a_p)_n}{(b_1)_n\cdots(b_q)_n}\frac{z^n}{n!},$$
with ${\bf a}=(a_1,\dots,a_p),$ ${\bf b}=(b_1,\dots,b_q)$ such that $-b_j\notin\mathbb{N}_0,$ $j\in\{1,\dots,q\},$ and $(a)_n=a(a+1)\cdots(a+n-1)={\Gamma(a+n)}/{\Gamma(a)}$ being the shifted factorial (or Pochhammer's symbol).
Our Theorem \ref{th2} on hyper-Bessel functions is another natural extension of Hurwitz' result and complements naturally the result of Chaggara and Ben Romdhane \cite[Theorem 4.2]{cha}. For simplicity we use the notation
\begin{align*}\mathcal{J}_{\mathbf{\alpha_d}}(z)=z^{-\frac{\alpha_1+\cdots+\alpha_d}{d+1}}J_{\mathbf{\alpha_d}}\left((d+1)\sqrt[d+1]{z}\right)=
\sum_{n\geq0}\frac{(-1)^nz^n}{n!\prod_{k=1}^d\Gamma(\alpha_k+n+1)}.\end{align*}

\begin{theorem}\label{th2}
All zeros of the hyper-Bessel function $\mathcal{J}_{\mathbf{\alpha_d}}$ are real when $\alpha_i>-1,$ $i\in\{1,\dots,d\}.$ Under the same conditions the function
$\mathcal{J}_{\mathbf{\alpha_d}}(-z)$ has only negative real zeros. Moreover, if $\alpha_i\geq0,$ $i\in\{1,\dots,d\},$ then $\mathcal{J}_{\mathbf{-\alpha_d}}(z)$ has exactly $[\alpha_1]+\cdots+[\alpha_d]$ non-positive zeros.
\end{theorem}

Recently, Kalmykov and Karp \cite{karp} conjectured that if $p<q,$ $\mathbf{b}>0,$ and $a_k>b_k$ for $k\in\{1,\dots,p\},$ then all zeros of the generalized hypergeometric function ${}_pF_q(\mathbf{a};\mathbf{b};z)$ are real and negative. Applying directly Laguerre's Lemma \ref{lemlag}, we conclude that
if $\mathbf{b}>0,$ then all zeros of ${}_0F_q(-;\mathbf{b};z)$ are real and negative. This particular result suggest the validity of the above
conjecture. Moreover, by a simple application of Lemma \ref{lemlag} it is possible to obtain the following result: if $p\leq q,$ $\mathbf{b}>0$ and $\mathbf{a}$ can be reindexed so that $a_k=b_k+n_k$ for $n_k\in\mathbb{N}$ and $k\in\{1,\dots,p\},$ then the function ${}_pF_q(\mathbf{a};\mathbf{b};z)$ has only negative real zeros. This result was stated by Richards \cite{ri}, who used Laguerre's Lemma \ref{lemlag}. The main idea is that
the function
$$\frac{\Gamma(a_1+z)\cdots\Gamma(a_p+z)}{\Gamma(b_1+z)\cdots\Gamma(b_q+z)}$$
is a meromorphic function as a quotient of two entire functions, however, after the reindexation the poles of the numerator are absorbed by those of the denominator, and hence it becomes an entire function of growth order $1$ for which Lemma \ref{lemlag} can be applied. We note that the reality of the zeros
of ${}_pF_q(\mathbf{a};\mathbf{b};z)$ follows also immediately from Obreschkoff's Lemma \ref{lemobr}. Moreover, by using P\'olya's Lemma \ref{lempolya} we obtain the following result, which complements \cite[Theorem 4]{karp}.

\begin{theorem}\label{th3}
Suppose that $\mathbf{b}>0$ and $\mathbf{a}$ can be reindexed so that $a_k=b_k+n_k$ for $n_k\in\mathbb{N}$ and $k\in\{1,\dots,p\}.$
\begin{enumerate}
\item[\bf 1.] If $p=q,$ then the function ${}_pF_q(\mathbf{a};-\mathbf{b};-z)$ has $[b_1]+\cdots+[b_p]+p$ non-positive zeros if $n_k\geq [b_k]+1$ for every $k\in\{1,\dots,p\};$ and has $n_1+\cdots+n_p$ non-positive roots if $n_k\leq[b_k]$ for each $k\in\{1,\dots,p\}.$
\item[\bf 2.] If $p<q,$ then the function ${}_pF_q(\mathbf{a};-\mathbf{b};-z)$ has $[b_1]+\cdots+[b_p]+[b_{p+1}]+\cdots+[b_q]+q$ non-positive zeros if $n_k\geq [b_k]+1$ for every $k\in\{1,\dots,p\};$ and has $n_1+\cdots+n_p+[b_{p+1}]+\cdots+[b_q]+q-p$ non-positive roots if $n_k\leq[b_k]$ for each $k\in\{1,\dots,p\}.$
\end{enumerate}
\end{theorem}

Now, we present some results for zeros of derivatives of Bessel functions. We note that the reality of the zeros stated in Theorem \ref{th4} was already proved in \cite{baricz} by using mathematical induction, and the rest of Theorem \ref{th4} is in agreement with \cite[Open Problem 1]{baricz}, which states that if $n-2s-2<\nu<n-2s-1,$ $s\in\mathbb{N}_0,$ then $J_{\nu}^{(n)}(z)$ has $4s+2$ complex zeros; while if $n-2s-1<\nu<n-2s,$ $s\in\mathbb{N},$ then $J_{\nu}^{(n)}(z)$ has $4s$ complex zeros.

\begin{theorem}\label{th4}
Let $n\in\mathbb{N}_0.$ If $\nu>n-1,$ then all zeros of $J_{\nu}^{(n)}(z)$ are real. Moreover, if $\nu\geq0,$ then $2^{\nu}z^{{\nu+n}}J_{-\nu}^{(n)}(2z)$ has at most $2[\nu]+2n$ complex zeros. In other words, if $n-2s-2<\nu<n-2s-1,$ $s\in\mathbb{N}_0,$ then $J_{\nu}^{(n)}(z)$ has at most $4s+2$ complex zeros; while if $n-2s-1<\nu<n-2s,$ $s\in\mathbb{N},$ then $J_{\nu}^{(n)}(z)$ has at most $4s$ complex zeros.
\end{theorem}

Now, we consider the functions $\Phi_\nu$ and $\Pi_\nu$, defined by
$$\Phi_\nu(z)=J_\nu(z)I_\nu'(z)-I_\nu(z)J_\nu'(z), \ \ \mbox{and}\ \ \Pi_\nu(z)=J_\nu(z)I_\nu(z),$$
where $I_\nu$ stands for the modified Bessel functions of the first kind. If $z\in \mathbb{C}$ and $\nu \in \mathbb{R}$ such that $\nu\neq -1,-2,\ldots$ then the functions $\Phi_\nu$ and $\Pi_\nu$ can
be written as follows (see \cite[p. 148]{watson}, \cite{abp}):
$$\Phi_\nu(z)=2\sum_{n\ge 0}\frac{(-1)^n(\frac{z}{2})^{2\nu+4n+1}}{n!\Gamma(\nu+n+1)\Gamma(\nu+2n+2)}$$
and
$$\Pi_\nu(z)=\sum_{n\ge 0}\frac{(-1)^n(\frac{z}{2})^{2\nu+4n}}{n!\Gamma(\nu+n+1)\Gamma(\nu+2n+1)}.$$
The following theorem on the functions
$$\mathcal{A}_{\nu}(z)= z^{\frac{-2\nu-1}{2}}\Phi_\nu(2\sqrt{z})
=2\sum_{n\ge 0}\frac{(-1)^nz^{2n}}{n!\Gamma(\nu+n+1)\Gamma(\nu+2n+2)}$$
and
$$\mathcal{B}_{\nu}(z)=z^{-\nu}\Pi_\nu(2\sqrt{z})
=\sum_{n\ge 0}\frac{(-1)^nz^{2n}}{n!\Gamma(\nu+n+1)\Gamma(\nu+2n+1)}$$
is another interesting application of Lemma \ref{lemlag} and \ref{lemobr} and is related to \cite[Open Problem 1]{bsy}.

\begin{theorem}\label{th5}
If $\nu>-1$, then all zeros of $\mathcal{A}_{\nu}(z)$ and $\mathcal{B}_{\nu}(z)$
are real. Moreover, if $\nu>-1$, then all the zeros of   $\mathcal{A}_{\nu}(\sqrt{z})$ and $\mathcal{B}_{\nu}(\sqrt{z})$ are real and positive. In addition, if $\nu\geq 0$, then $\mathcal{A}_{-\nu}(z)$ and $\mathcal{B}_{-\nu}(z)$ have at most $4[\nu]$ complex zeros.
\end{theorem}

We end this subsection with Theorem \ref{th6}. This result is related to \cite[Problem 6]{csord}: {\em is it true that for every $s\in\mathbb{R}^{+},$ there exists an $m\in\mathbb{N}$ such that $\sum\limits_{n\geq0}\frac{n^s}{(n!)^m}z^n\in\lp^{+}$?} Theorem \ref{th6} verifies the case when $s\in\mathbb{N}_0$ and $m=2.$

\begin{theorem}\label{th6}
If $s\in\mathbb{N}_0,$ then $\eta_{s}(z)=\sum\limits_{n\geq0}\frac{n^s}{\Gamma(n+1)}\frac{z^n}{n!}\in\lp^{+}.$
\end{theorem}

\subsection{Concluding remarks} We have seen that in some special cases we can find an upper bound for the exact number of complex zeros of entire functions. However, for example in the above mentioned conjecture of Kalmykov and Karp \cite[Conjecture 3]{karp} the coefficients are meromorphic and so far we are not aware on any result in the literature which would help to verify this conjecture. Craven and Csordas \cite[Problem 1.2]{craven} posed the following problem: {\em Characterize the meromorphic function $F$ with the property that $\sum\limits_{n\geq0}\frac{a_nF(n)}{n!}z^n$ is a transcendental entire function with only real zeros whenever the entire function $\sum\limits_{n\geq0}\frac{a_n}{n!}z^n$ has only real zeros.} This problem is strongly related to the above conjecture and its solution would give many new results on real and complex zeros of different special functions. It would be also of great interest to verify whether the condition on the difference of two consecutive zeros in Lemma \ref{lempolya} can be relaxed and to find the analogue of Lemma \ref{lempolya} when the function $G$ is meromorphic.

\section{Proofs of Main Results}
\setcounter{equation}{0}

\begin{proof}[\bf Proof of Theorem \ref{th1}]

First we consider the proof of the reality of the zeros of the Wright function $\phi(\rho,\beta,-z).$ The function $q_{\rho,\beta}:[0,\infty)\to\mathbb{R},$ defined by $q_{\rho,\beta}(z)=\frac{1}{\Gamma\left(\rho \frac{z}{2}+\beta\right)},$ is an entire function of order $1,$ belongs to $\mathcal{LP}$ and if $\rho,\beta>0$ then clearly it has no positive zero. By using Lemma \ref{lemobr} it follows that the function
$$\phi(\rho,\beta,-z^2)=\sum_{n\geq0}\frac{(-1)^nq_{\rho,\beta}(2n)}{n!}z^{2n}=\sum_{n\geq0}\frac{(-1)^nz^{2n}}{n!\Gamma(\rho n+\beta)}$$
has at most $0$ complex zeros, that is, all of its zeros are real. This implies that $\phi(\rho,\beta,-z)$ has also only real zeros when $\rho,\beta>0.$

An alternative proof of the fact that $\phi(\rho,\beta,-z)$ has only real zeros if $\rho,\beta>0$ is based on Runckel's Lemma \ref{lemrun}. Since $q_{\rho,\beta}(z)$ is of type \eqref{typeA} when $\rho,\beta>0,$ if we choose $f(z)=e^{-\left(\frac{z}{2}\right)^2},$ then by using Runckel's above mentioned result (that is, Lemma \ref{lemrun}) we obtain that the function $\phi(\rho,\beta,-z)$ has real zeros only if $\rho,\beta>0.$

Next, we show that if $0<\rho\leq1$ and $\beta>0,$ then all zeros of $\phi(\rho,\beta,-z)$ are positive. For this we consider the function $G_{\rho,\beta}:[0,\infty)\to\mathbb{R},$ defined by $G_{\rho,\beta}(z)=\frac{1}{\Gamma(\rho z+\beta)}.$ This function has zeros $z_k=-\rho^{-1}(\beta+k),$ where $k\in\mathbb{N}_0.$ These zeros are clearly simple and the distance between two consecutive zeros is
$\Delta_k=z_{k+1}-z_k=\rho^{-1}$ for every $k\in\mathbb{N}_0.$ The simplicity of the zeros is guaranteed by the Laguerre's theorem on separation of zeros (which states that if $f(z)$ is an entire function, not a constant, which is real for real $z$ and has only real zeros, and is of genus $0$ or $1,$ then the zeros of $f'$ are also real and are separated by the zeros of $f$) and by the fact the reciprocal gamma function is an entire function of genus $1.$ If $\rho\in(0,1],$ then $\Delta_k\geq1$ for every $k\in\mathbb{N}_0.$ Moreover, since $e^{-z}\in\mathcal{LPI}$ and the function $G_{\rho,\beta}$ has no zeros in $[0,\infty)$ when $\rho\in(0,1]$ and $\beta>0,$ by applying Lemma \ref{lempolya} we have that
$$\phi(\rho,\beta,-z)=\sum_{n\geq0}G_{\rho,\beta}(n)\frac{(-z)^n}{n!}=\sum_{n\geq0}\frac{(-1)^nz^n}{n!\Gamma(\rho n+\beta)}$$
has no non-positive real roots, that is, all its zeros are positive.

Now, we prove that the condition $\rho\leq 1$ can be relaxed. The growth order of the entire function $\phi(\rho,\beta,-z)$ is $(\rho+1)^{-1}$ (which is a non-integer number and lies in $(0,1)$ for $\rho>-1$) and thus in view of the Hadamard factorization theorem on growth order of entire functions it follows that for $\rho,\beta>0$ we have
$$\Gamma(\beta)\phi(\rho,\beta,-z^2)=\prod_{n\geq1}\left(1-\frac{z^2}{\lambda^2_{\rho,\beta,n}}\right),$$
where $\lambda_{\rho,\beta,n}$ denotes the $n$th positive zero of $\phi(\rho,\beta,-z^2),$ and this product is uniformly convergent on compact subsets of the complex plane. Consequently, we have that
$$\Gamma(\beta)\phi(\rho,\beta,-z)=\prod_{n\geq1}\left(1-\frac{z}{\lambda^2_{\rho,\beta,n}}\right),$$
which shows that all zeros of $\phi(\rho,\beta,-z)$ can be represented as squares, and thus indeed all zeros of $\phi(\rho,\beta,-z)$ are positive.

Alternatively, since for $\rho,\beta>0$ the zeros of $G_{\rho,\beta}$ are all negative, if we apply Lemma \ref{lemlag} we immediately obtain that for $\rho,\beta>0$ the Wright function $\phi(\rho,\beta,z)$ has only real and negative zeros. This means that for $\rho,\beta>0$ the Wright function $\phi(\rho,\beta,-z)$ has only real and positive zeros.

Finally, we proceed similarly as above, when we considered the case of $0<\rho\leq1$ and $\beta>0.$ If $0<\rho\leq1$ and $\beta<0$ then
exactly $[|\beta|]+1$ of the zeros $z_k=-\rho^{-1}(\beta+k),$ $k\in\mathbb{N}_0,$ are in $[0,\infty).$ Consequently, in view of Lemma \ref{lempolya} (with the function $e^{-z}\in\lpi$), the function $\phi(\rho,\beta,-z)$ has exactly $[|\beta|]+1$ non-positive zeros, that is, we proved that if $0<\rho\leq 1$ and $\beta<0,$ then $\phi(\rho,\beta,-z)$ has $[|\beta|]+1$ non-positive zeros. Now, changing $\beta$ to $-\beta$ we complete the proof.
\end{proof}

\begin{proof}[\bf Proof of Theorem \ref{th2}]
Consider the product of reciprocals of gamma functions appearing in Theorem \ref{th2}, that is,
$$\Theta_{\mathbf{\alpha_d}}(z)=\frac{1}{\Gamma(\alpha_1+z+1)\cdots\Gamma(\alpha_d+z+1)}.$$
In view of the representation
$$\frac{1}{\Gamma(z)}=ze^{\gamma z}\prod_{k\geq1}\left(1+\frac{z}{k}\right)e^{-\frac{z}{k}}$$ we have that
when $\alpha_i>-1,$ $i\in\{1,\dots,d\},$ the function $\Theta_{\mathbf{\alpha_d}}\left(\frac{z}{2}\right)$ is entire of growth
order $1$ and has no positive zeros. By using Lemma \ref{lemobr} we obtain that $\mathcal{J}_{\mathbf{\alpha_d}}(z^2)$ has at most $0$ complex zeros, that is,
all of its zeros are real. This implies that all zeros of $\mathcal{J}_{\mathbf{\alpha_d}}(z)$ are real when $\alpha_i>-1,$ $i\in\{1,\dots,d\}.$ Clearly, under the same conditions, the function $\Theta_{\mathbf{\alpha_d}}\left({z}\right)$ is also entire of growth
order $1$ and has no positive zeros, and applying Lemma \ref{lemlag} we obtain that $\mathcal{J}_{\mathbf{\alpha_d}}(-z)$ has only negative real zeros. Now, since for fixed $i\in\{1,\dots,d\}$ the reciprocal of $\Gamma(z-\alpha_i+1)$ has zeros $\xi_k=\alpha_i-1-k,$ $k\in\mathbb{N}_0,$ and $[\alpha_i]$ of these zeros are positive or zero, the distance between two arbitrary consecutive zeros is equal to $1,$ applying Lemma \ref{lempolya} (with the function $e^{-z}\in\mathcal{LPI}$) we obtain that
$$\mathcal{J}_{\mathbf{-\alpha_d}}(z)=\sum_{n\geq0}\Theta_{\mathbf{-\alpha_d}}(n)\frac{(-z)^n}{n!}$$ has exactly $[\alpha_1]+\cdots+[\alpha_d]$ non-positive zeros.
\end{proof}

\begin{proof}[\bf Proof of Theorem \ref{th3}]
First we suppose that $p=q.$ After reindexation the expression
$$\frac{\Gamma(a_1+z)\cdots\Gamma(a_p+z)}{\Gamma(b_1+z)\cdots\Gamma(b_q+z)}$$
will have zeros only as solutions of the equation
$$\prod_{k=1}^p\prod_{s=1}^{n_k}(b_k+z+s-1)=0,$$
that is, $\zeta_k=-b_k-s+1,$ where $s\in\{1,\dots,n_k\}$ and $k\in\{1,\dots,p\}.$ Now, if we replace $\mathbf{b}$ by $-\mathbf{b},$ then the above zeros clearly will change to $b_k-s+1,$ where $s\in\{1,\dots,n_k\}$ and $k\in\{1,\dots,p\}.$
If $n_k\geq [b_k]+1$ for $k\in\{1,\dots,p\}$ fixed, then we have that exactly $[b_k]+1$ zeros are positive; if $n_k=[b_k],$ then $[b_k]$ zeros are positive; and when $n_k\leq [b_k]-1,$ then $n_k$ number of zeros are positive. Thus, applying Lemma \ref{lempolya} (for $e^{-z}\in\mathcal{LPI}$) the result when $p=q$ follows.
The case when $p<q$ is similar to the case when $p=q.$ We just need to take care of the reciprocals of the remained expressions like $\Gamma(b_q+z).$
\end{proof}

\begin{proof}[\bf Proof of Theorem \ref{th4}]
Consider the entire function
$$2^{\nu}z^{\frac{n-\nu}{2}}J_{\nu}^{(n)}(2\sqrt{z})=\sum_{m\geq0}\frac{(-1)^m\Gamma(\nu+2m+1)z^m}{m!\Gamma(\nu+2m-n+1)\Gamma(\nu+m+1)}.$$
The function $$q_{\nu}(2z)=\frac{\Gamma(\nu+2z+1)}{\Gamma(\nu+2z-n+1)\Gamma(\nu+z+1)}$$ is entire since the poles of the numerator are absorbed by those of the denominator, and has growth order $1.$ The zeros of $q_{\nu}$ are of the form $\varsigma_k=\frac{k-1-\nu}{2},$ $k\in\{1,\dots,n\},$ and $\tau_s=-1-\nu-s,$ $s\in\mathbb{N}_0.$ If $\nu>n-1,$ $n\in\mathbb{N}_0,$ then clearly none of the above zeros is positive, and thus $q_{\nu}$ has $0$ positive zeros in this case. According to Obreschkoff's Lemma \ref{lemobr} it follows that $z^{n-\nu}J_{\nu}^{(n)}(2z)$ has at most $0$ complex zeros, that is, has only real zeros. Now, if we consider the entire function $q_{-\nu}(z),$ then for $\nu\geq0$ it has $[\nu]+n$ positive zeros. Thus, applying again Obreschkoff's Lemma \ref{lemobr} we conclude that
$$2^{\nu}z^{\nu+n}J_{-\nu}^{(n)}(2z)=\sum_{m\geq0}\frac{(-1)^mq_{-\nu}(2m)}{m!}z^{2m}$$
has at most $2[\nu]+2n$ complex zeros. In other words, if $\nu\leq0,$ then $2^{-\nu}z^{-\nu+n}J_{\nu}^{(n)}(2z)$ has at most
$2[-\nu]+2n$ complex zeros. This means that if $n-2s-2<\nu<n-2s-1,$ $s\in\mathbb{N}_0,$ then $J_{\nu}^{(n)}(z)$ has at most $4s+2$ complex zeros; while if $n-2s-1<\nu<n-2s,$ $s\in\mathbb{N},$ then $J_{\nu}^{(n)}(z)$ has at most $4s$ complex zeros.
\end{proof}

\begin{proof}[\bf Proof of Theorem \ref{th5}]
Consider the functions
$$
a_\nu(z)=\frac{1}{\Gamma(\nu+\frac{z}{2}+1)\Gamma(\nu+z+2)}\ \mbox{and}\ b_\nu(z)=\frac{1}{\Gamma(\nu+\frac{z}{2}+1)\Gamma(\nu+z+1)},
$$
which are entire of growth order $1$. Note that the zeros of $a_\nu(z)$ are of the form $z_k=2(-1-\nu-k)$ and $z_l=-2-\nu-l$, where $k, l \in \mathbb{N}_0$ and the zeros of $b_\nu(z)$ are of the form $z_s=2(-1-\nu-s)$ and $z_t=-1-\nu-t,$ where $s, t \in \mathbb{N}_0$. Therefore $a_\nu(z)$ and $b_\nu(z)$ have no positive zeros if $\nu>-1$. Now, appealing to Obreschkoff's Lemma \ref{lemobr}, we obtain that
$\mathcal{A}_{\nu}(z)$ and $\mathcal{B}_{\nu}(z)$ have at most $0$ complex zeros, that is, all of their zeros are real. Hence all the zeros of $\mathcal{A}_{\nu}(\sqrt{z})$ and $\mathcal{B}_{\nu}(\sqrt{z})$ are real and positive.

Alternatively, the positivity  of the zeros of $\mathcal{A}_{\nu}(\sqrt{z})$ and $\mathcal{B}_{\nu}(\sqrt{z})$  can be proved also by using Laguerre's Lemma \ref{lemlag}. For this we note that
$$a_\nu(2z)=\frac{1}{\Gamma(\nu+z+1)\Gamma(\nu+2z+2)}\ \mbox{and}\ b_\nu(2z)=\frac{1}{\Gamma(\nu+z+1)\Gamma(\nu+2z+1)}$$
are entire functions of order $1$ and they assume real values along the real axis and possess only negative zeros if
$\nu>-1$. Therefore in view of Laguerre's Lemma \ref{lemlag} we obtain that the entire functions
$$
u_\nu(z)=2\sum_{n\ge 0}\frac{z^{n}}{n!\Gamma(\nu+n+1)\Gamma(\nu+2n+2)}
$$
and
$$
v_\nu(z)=\sum_{n\ge 0}\frac{z^{n}}{n!\Gamma(\nu+n+1)\Gamma(\nu+2n+1)}
$$
have also real and negative zeros. Hence $u_\nu(-z)$ and $v_\nu(-z)$ have real and positive zeros. That is, $\mathcal{A}_{\nu}(\sqrt{z})$ and $\mathcal{B}_{\nu}(\sqrt{z})$ have real and positive zeros.

Now consider the entire functions $a_{-\nu}(2z)$ and $b_{-\nu}(2z)$. For $\nu\geq 0$,  $a_{-\nu}(2z)$ and $b_{-\nu}(2z)$ both have $2[\nu]$ positive zeros. Therefore by using the Obreschkoff's Lemma \ref{lemobr} we conclude that $\mathcal{A}_{-\nu}(z)$ and $\mathcal{B}_{-\nu}(z)$ have at most $4[\nu]$ complex zeros.
\end{proof}

\begin{proof}[\bf Proof of Theorem \ref{th6}]
Since the coefficients of $\eta_s$ are positive, we just need to show that $\eta_s(z)\in\lpi$ or equivalently $\eta_s(-z)\in\lpi.$ The entire function $z^s/\Gamma(z+1)$ has growth order $1,$ it belongs to $\mathcal{LP}$ and has only one zero in $[0,\infty),$ and that is $0.$ Consequently by using Lemma \ref{lempolya} (for the function $e^{-z}\in\lpi$), we obtain that $\eta_{s}(-z)$ has exactly $1$ nonpositive real root, and that is $0$. In other words, all zeros of $\eta_s(-z)$ are strictly positive, that is, $\eta_s(-z)\in\lpi.$
\end{proof}

\section*{Acknowledgement} The research of \'A. Baricz was supported by the STAR-UBB Advanced Fellowship-Intern of the Babe\c{s}-Bolyai University of Cluj-Napoca. Both authors wish to acknowledge the referee's comments and suggestions which enhanced this paper.

\end{document}